\documentclass[12pt]{article}
\usepackage{amsfonts,enumerate,amsmath}

\newtheorem{theorem}{Theorem}

\newtheorem{proposition}{Proposition}

\def\R{{\mathbb R}}
\def\Z{{\mathbb Z}}

\def\beq{\begin{equation}}
\def\eeq{\end{equation}}

\def\sn{{\mathrm{sn}\,}}

\begin{document}

\title{On an integrable magnetic geodesic flow on the two-torus
\thanks{The work was supported by RSF (grant 14-11-00441).}
}
\author{I.A. Taimanov
\thanks{Sobolev Institute of Mathematics, 630090 Novosibirsk, Russia, and
Department of Mechanics and Mathematics, Novosibirsk State University,
630090 Novosibirsk, Russia
e-mail: taimanov@math.nsc.ru.} }
\date{}
\maketitle

\section{Introduction}

The motion of a charge in a magnetic field on a configuration space $M^n$
is described by the Euler--Lagrange equations for the Lagrangian function
$$
L(q,\dot{q}) = \frac{1}{2}g_{ik}\dot{q}^i \dot{q}^k + A_i \dot{q}^i
$$
where the first term is the kynetic energy calculated by using the Riemannian metric $g_{ik}$ on $M^n$
and the second term describes the interaction of a charge with the magnetic field $F$ which is a closed two-form $F$ on
$M^n$ such that $F=dA$. 
If $F$ is non-exact, then the one-form $A$ is defined only locally.
If $F=0$, then we get the Lagrangian function for geodesics on $M^n$ and for this reason solutions $q(t)$
for the Euler--Lagrange equations with a general magnetic field are called magnetic geodesics.

Trajectories lying on the energy level 
$$
E = \frac{1}{2}g_{ik}\dot{q}^i \dot{q}^k = \mathrm{const}
$$
satisfy to the Euler--Lagrange equations for the reduced Lagrangian function 
$$
L_E(q,\dot{q}) = \sqrt{2E} \sqrt{g_{ik}\dot{q}^i \dot{q}^k} + A_i \dot{q}^i.
$$

In this article we consider the magnetic geodesic flow on the two-torus
$$
T^2 = \R^2/(2\pi \Z)^2,
$$
where the coordinates $x$ and $y$ are defined modulo $2\pi$, with the flat Riemannian metric
$$
ds^2 = dx^2 + dy^2
$$
and with the magnetic field given by the $2$-form
$$
F = \cos x \, dx \wedge dy.
$$
This form is exact:
$$
F = dA \ \ \ \mbox{with}\ \ \ A = \sin x \, dy.
$$
Therewith we

\begin{enumerate}
\item
{\sl completely integrate the flow and, in particular, describe all trajectories in terms of elliptic functions} (Theorem \ref{thint});

\item
{\sl show that for all contractible periodic magnetic geodesics the reduced action
$$
S_E(\gamma)  = \int_\gamma L_E \, dt
$$
is positive:}
$$
S_E > 0;
$$
(Theorem \ref{thcontr});

\item
{\sl for $E<\frac{1}{2}$ find the minimizers of the action $S_E$ extended to submanifolds $\Sigma$ of $T^2$ as follows:}
$$
S_E(\Sigma,f) = \sqrt{2E}\, \mathrm{length}\,(\partial f(\Sigma)) + \int_\Sigma f^\ast (F),
$$
where $f: \Sigma \to T^2$ is the embedding
(Theorem \ref{minfilm}
\footnote{This extension of $S_E$ for the space of films $S_E$ was introduced in \cite{T1,T2}. Theorem \ref{minfilm}
describes the minimal films for the case $E<\frac{1}{2}$, i.e. when $S_E$ attains negative values.});

\item
{\sl explicitly describe all contractible periodic magnetic geodesics and, in particular, show that they exist only for
$E < \frac{1}{2}$ and simple (not iterates of others) contractible periodic magnetic geodesics form two $S^1$-families}
(Theorem \ref{contractible}).
\end{enumerate}

The initial intention for writing this article was to supply the study of periodic problem for magnetic geodesics with a non-trivial explicit example of a magnetic geodesic flow
with interesting properties. 

The study of the periodic problem for magnetic geodesics was initiated by Novikov in early 1980s \cite{NS,Novikov1981,Novikov} in the framework of qualitative study of certain Hamiltonian systems from classical mechanics. It appears that an application of the classical Morse theory to proving the existence of periodic magnetic geodesics
meets many obstacles which can not be overgone by classical methods \cite{Novikov}.
We gave a survey of that as well as of the first results in this direction in our survey \cite{T0}.

For dealing with these difficulties new ideas methods arising from symplectic geometry, dynamical systems and fixed points theory were proposed
(see, for instance, \cite{A,BT,Benedetti,Branding,CMP,Contreras,FS,Ginzburg,Ginzburg2,GG,Kozlov,NT,Schneider,T83,T1,T2}).

One of them consists in a proposal to extend the reduced action functional $S_E$ for closed curves in an oriented two-manifold to the space of films, and then to establish the existence of a minimal film whose boundary components would be locally minimal periodic magnetic geodesics. That was done in \cite{T1,T2} where the existence of nontrivial periodic magnetic geodesics was established for strong magnetic fields. 

Let us assume that a magnetic field $F$ is exact. Then it is easy to satisfy the strength condition from \cite{T1,T2}
by multiplying $F$ by a sufficiently large constant $\lambda$ or by considering only low energy levels $E$. In fact, the ratio $\frac{\lambda}{\sqrt{E}}$ has to be sufficiently large, or, in other words, $E <E_0$ where $E_0$ is some constant. However if $E>E_1$ with $E_1$ a constant, then $S_E$ is the length functional for a Finsler metric and the periodic problem can be studied by the classical Morse theory. It was established in \cite{CMP} that $E_0=E_1=C$, the constant $C$ is the Mane critical level and on this level the
existence of a periodic magnetic geodesic is derived from the Aubry--Mather theory. Hence the existence of a periodic magnetic geodesic was established for all energy levels 
for exact magnetic fields on oriented closed two-manifolds. 

Recently the existence of infinitely many periodic magnetic geodesics was established for exact magnetic fields on two-manifolds  for almost all subcritical energy levels, i.e. for almost all $E<C$ \cite{AMP,AMMP}. For other recent results which clarify the periodic problem on two-manifolds we refer to \cite{AB,BL,BZ}.

To finish the introduction we are left to make two

{\sc Remarks.} 1) Given an exact magnetic field, if the functional $S_E$ attains negative values on certain closed contractible curves, then the whole manifold of one-point curves
$M^n$ can be overthrown into domain $\{S_E <0\}$ and hence every $k$-cycle from $M^n$ generates a nontrivial cycle from
$H_{k+1}(P_0(M^n),\{S_E \leq 0\})$ where $P_0$ is the space of contractible closed curves in $M^n$ ({\it ``the principle of throwing out cycles''}) \cite{T83,T10}.
It would be interesting to understand how the contractible periodic magnetic geodesics given by Theorem \ref{contractible} fit into this picture and 
how to compute $H_\ast(P_0(M^n),\{S_E \leq 0\})$ by using them.

2) The existence of a minimal film proved in \cite{T1,T2} does not say anything about the existence of contractible periodic magnetic geodesics.
A priori if we have a nonselfintersecting closed curve $\gamma$ with $S_E(\gamma)<0$, then it is not necessary that there exists a minimizer of $S_E$ in the class
of nonselfintersecting closed curves with $S_E<0$. Theorems \ref{thcontr} and \ref{minfilm} say that the flow in study supplies an explicit example of such a situation.
Actually this article arises as an answer to the question by A. Abbondandolo on such an example asked by him after my talk in Bochum in May 2015.

\section{The flow and the first integrals}

The motion of a charge in $T^2$ in the magnetic field $F$ is defined by the Lagrangian function
$$
L(q,\dot{q}) = \frac{1}{2} g_{ik}\dot{q}^i\dot{q}^k - A_i\dot{q}^i = \frac{\dot{x}^2+\dot{y}^2}{2} + \sin x \, \dot{y},
$$
where $q = (x,y) \in T^2$.

The corresponding Lagrangian system takes the form
\begin{equation}
\label{flow}
\begin{split}
\ddot{x} = \cos x\, \dot{y}, \\
\ddot{y} = -\cos x\, \dot{x}.
\end{split}
\end{equation}

Since the Lagrangian function is independent on $y$, we have

\begin{proposition}
The system (\ref{flow}) has two functionally independent first integrals:

1) the energy
$$
E = \frac{\dot{x}^2+\dot{y}^2}{2},
$$

2) the momentum corresponding to $y$:
$$
p:= p_y = \dot{y} + \sin x.
$$
\end{proposition}

By this Proposition, we have
$$
\dot{y} = p - \sin x, \ \ \ \dot{x}= \sqrt{2E - (p-\sin x)^2},
$$
and therefore these equations are integrable in quadratures:
\begin{equation}
\label{quadrature}
\begin{split}
t  = \int_{x_0}^x \frac{dv}{\sqrt{2E - (p-\sin v)^2}}, \\
y(t) = y_0 + pt - \int_{0}^t \sin x(\tau)\,d\tau,
\end{split}
\end{equation}
where the constants $p$ and $E$ are defined by the initial data:
$$
x_0 = x(0), \ \ y_0 = y(0), \ \ \dot{x}(0), \ \ \dot{y}(0).
$$
To find $x(t)$ we have to invert the the integral in the first equation from (\ref{quadrature}.
To do that we make the substitution
$$
z = \sin x
$$
and derive
$$
t = \int_{\sin x_0}^{\sin x} \frac{dz}{w}
$$
where
\begin{equation}
\label{elliptic}
w^2 = (1-z^2)(2E - (p-z)^2).
\end{equation}
This elliptic curve is reduced to the Legendre normal form by the standard procedure (see \S \ref{section_Legendre})
and $x(t)$ is expressed in terms of the Jacobi elliptic functions.

\begin{theorem}
\label{thint}
$$
\sin x(t) = \sn \left(\frac{t+D}{C},k\right),
$$
$$
y(t) = y_0 + pt - \int_{0}^t \sn \left(\frac{\tau+D}{C},k\right) \,d\tau,
$$
where
$D = \int_0^{\sin x_0} \frac{dz}{w}$ and the constants $k$ and $D$  are determined by the reduction of the elliptic curve (\ref{elliptic}) to the
Legendre normal form (see (\ref{constk}) and (\ref{constC})).
\end{theorem}

{\sc Remark.} This flow is a simplest example of a magnetic analog of the geodesic flow on a surface of revolution. The general family 
is given by pairs consistsing of a metric $g(x)(dx^2 + dy^2)$ and a two-form $h^\prime(x)\, dx\wedge dy$ where the functions $g$ and $h$ are both periodic.
In this case the additional conservation law needed for complete integrability takes the form of the ``Clairaut integral'':
$$
\frac{\partial L}{\partial \dot{y}} = \frac{\dot{y}}{g(x)} + h(x).
$$
We think that in this class one can find other interesting examples of magnetic geodesic flows.

\section{The variational principle for closed trajectories}
\label{section_Var}

The Lagrangian is not homogeneous in velocities and therefore
the restrictions of the system on different energy levels are not trajectorically isomorphic
as in the case of the geodesic flow. In fact, the restriction of the system on the energy level is described by the Lagrangian
$$
L_E = \sqrt{2E}\sqrt{g_{ik}\dot{q}^i\dot{q}^k} + A_i \dot{q}^i = \sqrt{2E}\sqrt{\dot{x}^2+\dot{y}^2} + \sin x \, \dot{y}.
$$
Indeed, the Lagrangian $L_E$ is homogeneous of first order in velocities and the corresponding trajectories are defined up
to parameterizations. The Euler--Lagrange equations takes the form
$$
\frac{d}{dt}\left(\frac{\sqrt{2E}g_{ik}\dot{q}^k}{\sqrt{g_{ik}\dot{q}^i\dot{q}^k}}\right) = \frac{\partial L_E}{\partial q^i}
$$
and for $E = \frac{1}{2}g_{ik}\dot{q}^i\dot{q}^k$ both left- and -right-hand sides reduces to the left- and right-hand sides for the
Euler--Lagrange equation for $L$.

The following statement is easily checked by straightforward computations

\begin{proposition}
\label{mane}
The Lagrangian function $L_E(q,\dot{q})$ is positive for all $\dot{q}\ne 0$ if and only if $E > \frac{1}{2}$.

If $E=\frac{1}{2}$, then $L_E$ vanishes at $\{\dot{x}=0, \dot{y}>0, \sin x =  -1\}$ and $\{\dot{x}=0, \dot{y}<0, \sin x = 1\}$.

$L_E$ takes negative values for every $E < \frac{1}{2}$.
\end{proposition}

{\sc Remark.} The particular level $E = \frac{1}{2}$ may be treated as the Mane critical level of the system. For magnetic geodesic flows
the Mane critical level $C$ is defined for magnetic fields $F$ on Riemannian manifolds $M$ such that the pullback $\pi^\ast F$
of $F$ onto the universal covering $\widetilde{M}$ of $M$ is exact.  In this case it is defined as
$$
C = \inf_{\theta} \sup_{{M}} H(q,\theta)
$$
where $H(q,p) = \frac{1}{2}|p|^2$ is the Hamiltonian function of the system and $d\theta = \pi^\ast F$. In our case we have
$$
C = \inf_{\theta} \sup_{\R^2} \frac{1}{2}(\theta_1^2 + \theta_2^2)
$$
where $\theta = \theta_1 \, dx + \theta_2 \, dy$. Since $d\theta = \cos x \, dx \wedge dy$, we have
$\theta_1 = f_x, \theta_2 = f_y+\sin x$, where $df = f_x \, dx+ f_y\, dy$ and $f: \R^2 \to \R$ is a smooth function.
The restriction of $f_y$ onto the line $\sin x = 1$ is $2\pi$-periodic and $\int_0^{2\pi}f_y dy =0$. Therefore on this line $f_y$ achieves its minimum at which $f_y=0$ and
$(f_y+\sin x)^2=1$ . That implies that $\sup_{\R^2} H(q,\theta) \geq \frac{1}{2}$. This lower bound becomes exact at $f \equiv 0$ and therefore
$C = \frac{1}{2}$.

Closed trajectories of the flow lying on the energy level $E$ are extremals of the functional
$$
S_E(\gamma) = \int_\gamma L_E\, dt
$$
defined on the space of closed curves. However this functional is not always bounded from below in contrast to the length functional studied in the classical Morse theory
and moreover for non-exact magnetic fields, i.e. when the closed $2$-form is not globally  represented as a differential $dA$ of a $1$-form, this functional
is multi-valued. These differences with the Morse they were discussed in detail by Novikov \cite{Novikov} (see also \cite{T0}).

In our case the magnetic field $F$  is exact however, by Proposition \ref{mane}, $S_E$ is not bounded from below for small values of $E$.

The variational study of the periodic problem is hindered also by the fact that fit is not known does $S_E$ satisfies the Palais--Smale type conditions.
 A priori the deformation decreasing $S_E$ may diverge even when $S_E$ approaches a critical level.

In \cite{T1,T2}  for studying periodic magnetic geodesics in two-dimensional manifolds
we introduced the spaces of films which are embeddings $f: \Sigma^2 \to M^2$ of oriented two-manifolds $\Sigma^2$
 with boundaries into a two-dimensional closed
 manifold $M^2$ endowed with a Riemannian metric and with a magnetic field $F$ which is not necessarily exact.
 For such films
 $$
 f: \Sigma^2 \to M^2
 $$
 there is defined an action functional
 $$
 S_E(\Sigma,f)  = \sqrt{2E}\, \mathrm{length}\, (f(\partial \Sigma))  + \int_{\Sigma^2}f^\ast (F),
 $$
 which for exact magnetic fields reduces to $S_E = \int_{f(\partial \Sigma^2)} L_E \, dt$.
 It was proved in \cite{T1,T2} that

 \medskip

 {\sl if
 \begin{enumerate}
 \item
 $F$ is exact,  or
 $F$ is non-exact and $\int_{M^2} F `>0$ (as we assume without loss of generality),

 \item
 there is a film with $S_E <0$,
 \footnote{Such a field is called {\it strong} if $F$ is exact or {\it oscillating} if $F$ is non-exact.}
 \end{enumerate}

 \noindent
 then there exists a film $\Sigma$
 on which $S_E$ attains its minimal value among films and the boundary of
 $\Sigma$ consists of closed magnetic trajectories which are locally minimal for $S_E$}.

\section{Contractible periodic magnetic geodesics}

By (\ref{quadrature}),
$$
y(t) = y(0) + \int \frac{(p-z)\,dz}{\sqrt{(1-z^2)(2E - (p-z)^2)}}, \ \ \ z=\sin x(t).
$$
Therefore for the closed trajectory $\gamma(t)$ of the flow we have
$$
\Delta y = \int_\gamma \dot{y}\,dt = \int_\gamma \frac{(p-z)\,dz}{\sqrt{(1-z^2)(2E - ((p-z)^2)}}
$$
where $\Delta y$ is the increment of $y$ along the pullback of the trajectory onto the universal covering $\R^2 \to T^2$.
If a closed trajectory is contractible, then
$$
\Delta y =0.
$$
Let us  compute
$$
S_E(\gamma) = \int_\gamma L_E dt =
\int_\gamma (\sqrt{2E}\sqrt{\dot{x}^2 + \dot{y}^2} + \sin x \, \dot{y})\,dt =
$$
$$
= \int_\gamma (\sqrt{2E}\sqrt{\dot{x}^2 + \dot{y}^2} + \sin x \, \dot{y})\,dt =
 \int_\gamma (\sqrt{2E}\sqrt{\dot{x}^2 + \dot{y}^2} + (p - \dot{y})\dot{y})\,dt =
 $$
 $$
= \int_\gamma (\sqrt{2E}\sqrt{\dot{x}^2 + \dot{y}^2} - \dot{y}^2))\,dt  + p \int_\gamma \dot{y}\,dt.
$$
Since $\sqrt{2E}\sqrt{\dot{x}^2 + \dot{y}^2} = \dot{x}^2 + \dot{y}^2$, $p$ is constant, and $\int_\gamma \dot{y}\,dt = \Delta y$,
we have
\begin{equation}
\label{increment}
S_E(\gamma) = \int_\gamma \dot{x}^2\, dt  + p \Delta y \geq 0.
\end{equation}

\begin{proposition}
Given the energy level $E$, if $\gamma$ is a nontrivial (different from a one-point contour) closed magnetic geodesic with $\dot{x} \equiv 0$,
then it is one of the following orbits:
$$
\gamma_{\pm \pm} = \left\{ x(t) = \pm \frac{\pi}{2}, \ y(t) = \pm \sqrt{2E}\, t, \ 0 \leq t \leq \pi \sqrt{\frac{2}{E}}\right\},
$$
and no one of these orbits is contractible.
\end{proposition}

The proof of this Proposition immediately follows from (\ref{flow}).
Together with (\ref{increment}) this Proposition implies

\begin{theorem}
\label{thcontr}
For a contractible periodic magnetic geodesic $\gamma$ which is different from a one-point contour and lies on the energy level $E$,
we have
$$
S_E(\gamma)  = \int_\gamma \dot{x}^2\,dt > 0.
$$
\end{theorem}

Let us consider a film $\Pi$ formed by the embedding of the cylinder
$$
\Pi = \left\{ \frac{\pi}{2} \leq x \leq \frac{3\pi}{2}, \ \ 0 \leq y \leq 2\pi\right\}
$$
into $T^2$. Its image is the closure of the domain on which $F < 0$.
The boundary of $\Pi$ is formed by a pair of closed trajectories:
$$
\partial \Pi = \gamma_{-+} \cup \gamma_{+-}
$$
and it is easy to check that
$$
S_E (\Pi) = \sqrt{2E}\, \mathrm{length} \, (\partial \Pi) + \int_0^{2\pi} dy \int_{\pi/2}^{3\pi/2} \cos x dx =
4\pi(\sqrt{2E}-1).
$$

\begin{theorem}
\label{minfilm}
For every $E < \frac{1}{2}$ the functional $S_E$ on the space of films attains its minimal value on $\Pi$.
\end{theorem}

{\sc Proof.}
By the results of \cite{T1} mentioned in \S \ref{section_Var}, $S_E$ attains its minimum on a film $\Sigma$
whose boundary components are local minima of $S_E$.

Let all the boundary components are contractible. This is possible if $\Sigma$ consists in components diffeomorphic to discs with holes.
The boundary of every hole is a contractible magnetic geodesic $\gamma$ and, by Theorem \ref{thcontr}, $S_E >0$ on such a contour.
If we glue the hole by a disc, we obtain a new film $\widetilde{\Sigma}$ with
$$
S_E (\widetilde{\Sigma} ) = S_E(\Sigma) - S_E(\gamma),
$$
i.e. decrease the value of $S_E$ which contradicts to the definition of $\Sigma$ as a global minimum.

Hence $\partial \Sigma$ contains a non-contractible component $\gamma_1$ and, since $\partial \Sigma$ realizes a trivial class in $1$-homologies,
it has to contain another non-contractible component $\gamma_2$. Therefore
$$
S_E(\Sigma) \geq \mathrm{length}\,(\gamma_1) + \mathrm{length}\, (\gamma_2) + \int_\Psi F.
$$
However the lengths of noncontractible contours are at least $\sqrt{2E}\,2\pi$, i.e., the length of $\gamma_{+-}$ and of $\gamma_{-+}$
which are minimal non-contractible geodesics on $T^2$.
Therefore
$$
S_E(\Sigma) \geq 4\pi + \int_\Sigma F
$$
and the right-hand side achieves its minimum on the film $\Pi$. Q.E.D.

Now let us describe all nontrivial (different from a point) contractible closed trajectories.

By (\ref{quadrature}), every periodic in $x$ trajectory
is obtained by the inversion of the integral
\begin{equation}
\label{inversion}
\int \frac{dz}{\sqrt{(1-z^2)(2E-(p-z)^2)}}, \ \ \ z = \sin x,
\end{equation}
where $z$ goes along the bounded real oval, i.e. the contour, on the Riemannian surface
$$
w^2 = P(z) = (1-z^2)(2E-(p-z)^2),
$$
which covers the interval $I$ such that
$P(z) \geq 0$ on the interval and $P(z)$ vanishes at its ends.
It is clear that
$$
I \subset [-1,1].
$$
Moreover it is clear that all roots of $P(z)$ are different otherwise the integral (\ref{inversion}) diverges and
does not correspond to an $x$-periodic solution.

On every periodic trajectory there exist a pair of points $q_1, q_2$ such that
\begin{equation}
\label{pair}
\dot{x}(q_1) = \dot{x}(q_2) = 0, \ \ \dot{y}(q_1)>0, \ \ \dot{y}(q_2)<0.
\end{equation}
The condition $\dot{x} = 0$ is equivalent to
\begin{equation}
\label{roots}
2E - (p- z)^2 = 0.
\end{equation}
If $I$ contains only one root of this equation then at the corresponding points on
an $x$-periodic trajectory
$\dot{y} = p-z$ has the same values and therefore there are no points $q_1$ and $q_2$
meeting (\ref{pair}). In this case the $x$-periodic trajectory does not close up and there is a nontrivial
translation period in $y$.
Hence for a periodic trajectory the roots $z_1 < z_2$ of (\ref{roots}) lie inside $I$:
$$
-1 < z_1 < z_2 < 1.
$$

We are left to check the last condition that the translation period in $y$ vanishes.
By (\ref{quadrature}), it is equal to
$$
\Delta y = 2 \int_{z_1}^{z_2} \frac{(p-z) \, dz}{\sqrt{(1-z^2)(2E-(p-z)^2)}} =
$$
$$
= 2 \int_{z_1}^{z_2} \frac{(p-z) \, dz}{\sqrt{(1-z^2)(z-z_1)(z_2-z)}}, \  \ \ z_1+z_2 = 2p.
$$
We have
$$
z_1 = p-a, \ \ \ z_2 = p+a, \ \ a>0,
$$
and by substitution $u=p-z$ we
derive that
$$
\Delta y = 2 \int_{-a}^a \frac{u\,du}{\sqrt{(1-(p-u)^2)(a^2-u^2)}}.
$$
By comparing the values of the integrand at $\pm u$, we infer that
$$
\Delta y \ \ \
\begin{cases} >0 & \text{for $p<0$}\\ <0 & \text{for $p>0$} \\ =0 & \text{for $p=0$}.
\end{cases}
$$

Let us summarize these facts in the following

\begin{theorem}
\label{contractible}
\begin{enumerate}
\item
For $E \geq \frac{1}{2}$ there are no contractible closed orbits.

\item
For every $E$ such that $0 < E < \frac{1}{2}$, there exist two $S^1$-families of
simple periodic magnetic geodesics. These families are invariant with respect to translations by $x$: $x \to x + \mathrm{const}$,
and obtained by the inversion of the integral
$$
t = \int \frac{dz}{\sqrt{(1-z^2)(2E-z^2)}}, \ \ \  -\sqrt{2E} \leq z = \sin x(t) \leq \sqrt{2E},
$$
and by solving the equation for dynamics in $y$:
$$
\dot{y} = - \sin x.
$$
All other nontrivial contractible periodic magnetic geodesics are iterates of these simple closed magnetic geodesics.

\item
These families lie in the domains separated by the contours on which $\sin x = \pm 1$.
In particular, no one contractible closed orbit intersects these contours.

\item
These families degenerate to the pair of contours $\{ \sin x = 0 \}$ formed by one-point closed curves as $E \to 0$.

\item
For these simple periodic magnetic geodesics
$$
S_E = 2 \int_{-a}^a \sqrt{2E-\sin^2 x}\, dx, \ \ \ a = \arcsin \, \sqrt{2E}.
$$
\end{enumerate}
\end{theorem}

Statement 4 is quite evident from the physical point of view: for very small energies closed orbits are trapped near
critical points of the magnetic field, i.e. of the function $f$ such that $F = f dx \wedge dy$.

Statement 5 is derived by straightforward computations from Theorem \ref{thcontr}.

\section{Appendix: The Legendre normal form of an elliptic curve and elliptic integrals}
\label{section_Legendre}

In this section we recall some facts on the reduction of an elliptic curve to the Legendre normal form
(for more details see, for instance, \cite{BE}). This is necessary for deriving Theorem \ref{thint}.

Let
$$
P(z) = (z-a_1)(z-a_2)(z-a_3)(z-a_4),
$$
be a polynomial with four different real zeroes $a_1,\dots,a_4$.
We recall how to transform the Riemann surface (elliptic curve)
\begin{equation}
\label{curve}
w^2 = P(z)
\end{equation}
to the Legendre form
\begin{equation}
\label{legendre}
\eta^2 = (1-\xi^2)(1-k^2 \xi^2).
\end{equation}

We enumerate the zeroes as follows:
$$
a_3 < a_1 < a_2 < a_4
$$
and decompose $P(x)$ into a product $P(z)=Q_1(z)Q_2(z)$ of two quadratic polynomials of the form
$$
Q_1(z) = (z-a_1)(z-a_2), \ \ \ Q_2(z) = (z-a_3)(z-a_4).
$$

Let us consider two cases:

1)  $Q_1(z) = z^2 -a_1^2$, $Q_2(z) = z^2-a_3^2$. Then the transformation
$$
\xi = \frac{z}{a_1}, \ \ \ \eta = \frac{w}{a_1^2 a_3^2}
$$
reduces the equation (\ref{curve}) to the form (\ref{legendre}) with $k^2 = \frac{a_1^2}{a_3^2}$.

2) If the case 1) does not hold,  then there exist $\lambda_1$ and $\lambda_2$ such that
$$
Q_1(z) - \lambda_1 Q_2(z) = (1-\lambda_1)(z-\mu)^2, \ \ \ \
Q_1(z) - \lambda_2 Q_2(z)  = (1-\lambda_2)(z-\nu)^2.
$$
These constants $\lambda_{1,2}$ are determined as the eigenvalues of the pair of quadratic forms defined by $Q_1$ and $Q_2$, i.e.,
as the zeroes of the equation
$$
\det \left(\begin{array}{cc} 1 - \lambda & -\frac{a_1+a_2}{2} + \lambda \frac{a_3+a_4}{2} \\
 -\frac{a_1+a_2}{2} + \lambda \frac{a_3+a_4}{2} & a_1 a_2 - \lambda a_3 a_4
 \end{array}\right) = 0.
 $$
 Therewith we have
 \begin{equation}
 \label{qsum}
 Q_1(z) = B_1 (z-\mu)^2 + C_1 (z-\nu)^2, \ \ \ Q_2(z) = B_2 (z-\mu)^2 + C_2 (z-\nu)^2.
 \end{equation}
 The constants $B_1,B_2,C_1,C_2$ are trivially computed as
 \begin{equation}
 \label{bc}
 B_j = \frac{Q_j(\nu)}{(\mu-\nu)^2}, \ \ \ C_j = \frac{Q_j(\mu)}{(\mu-\nu)^2}, \ \ \ j=1,2,
 \end{equation}
 and a substitution of  them into the equations $B_1+C_1=B_2+C_2=1$ leads to the following relations
 \begin{equation}
 \label{harmonicity0}
 2(\mu \nu + a_1 a_2) = (\mu + \nu)(a_1 + a_2), \ \ \
 2(\mu \nu + a_3 a_4) = (\mu + \nu)(a_3 + a_4).
 \end{equation}
 Since $\mu$ and $\nu$ are different from the zeroes $a_j, j=1,\dots,4$, the latter relations are rewritten as
 \begin{equation}
 \label{harmonicity}
 \frac{\nu - a_1}{\nu-a_2} = - \frac{\mu -a_1}{\mu-a_2}, \ \ \
 \frac{\nu - a_3}{\nu-a_4} = - \frac{\mu-a_3}{\mu-a_4}.
 \end{equation}
 The solutions $\mu$ and $\nu$ to (\ref{harmonicity})  are obtained as the zeroes of the quadratic equation:
 $$
 \lambda^2 - A\lambda + B = 0
 $$
 where, by (\ref{harmonicity0}), we have
 $$
 A = 2 \frac{a_1 a_2 - a_3 a_4}{a_1 + a_2 - a_3 -a_4}, \ \ \
 B = \frac{a_1 a_2 (a_3 + a_4) - a_3 a_4 (a_1 + a_2)}{a_1 + a_2 - a_3 - a_4}.
 $$
 It is easy to check that  $\mu$ and $\nu$
 have to correspond to  different real ovals of (\ref{curve}). This means, without loss of generality, that
 $$
 a_1 < \nu < a_2
 $$
 and
 $$
 \mu< a_3 \ \ \ \ \mbox{or} \ \ \ \ \mu> a_4.
 $$
 By (\ref{qsum}), we have
 $$
 w^2 = (B_1 (z-\mu)^2 + C_1 (z-\nu)^2)(B_2 (z-\mu)^2 + C_2(z-\nu)^2) =
 $$
 $$
 = B_1 (z-\mu)^2 \left(1 + \frac{C_1}{B_1} \frac{(z-\nu)^2}{(z-\mu)^2}\right)
 B_2 (z-\mu)^2 \left(1 + \frac{C_2}{B_2} \frac{(z-\nu)^2}{(z-\mu)^2}\right) =
 $$
 $$
 = B_1 B_2 (z-\mu)^4 (1-\xi^2)(1 - k^2 \xi^2)
 $$
 for
 $$
 \xi = \sqrt{-\frac{C_1}{B_1}} \frac{z-\nu}{z-\mu}, \ \ \ \ k^2 = \frac{B_1 C_2}{B_2 C_1}.
 $$
 We are left to put
 $$
 \eta = \frac{w}{\sqrt{B_1 B_2} (z-\mu)^2}
 $$
 to reduce the curve (\ref{curve}) to the Legendre normal form (\ref{legendre}).
 By (\ref{bc}) and (\ref{harmonicity}), we have
 $$
 \sqrt{B_1B_2} = \frac{\sqrt{P(\nu)}}{(\mu-\nu)^2} > 0,
 $$
 \begin{equation}
 \label{constk}
 k^2 = \left(\frac{\nu-a_1}{\mu-a_1}\right)^2 \left(\frac{\mu-a_3}{\nu-a_3}\right)^2,
 \end{equation}
 $$
 \lambda = \sqrt{-\frac{B_1}{C_1}} = \sqrt{\frac{(\nu-a_1)(\nu-a_2)}{(\mu-a_1)(a_2-\mu)}} ,
 $$
 $$
 \frac{1}{\lambda} \frac{a_1-\nu}{a_1-\mu} = \xi(a_1) = \pm 1.
 $$
 Since $\lambda$ is defined up to a sign, let us put
 $$
 \lambda = \frac{a_1-\nu}{\mu-a_1}
 $$
 to achieve
 $$
 \xi(a_1) = -1.
 $$
 Therewith we have
 $$
 \xi =  \frac{1}{\lambda} \frac{z-\nu}{z-\mu} = \frac{\mu-a_1}{a_1-\nu} \frac{z-\nu}{z-\mu},
 $$
 $$
 \xi(a_2) = 1, \ \ \ \xi(a_3) = -\frac{1}{k}, \ \ \ \xi(a_4) = \frac{1}{k} \ \ \ \mbox{where $k >0$,}
 $$
 and, since the real ovals $\{(z,w), P(z) \geq 0\}$ are mapped into real ovals by the transformation $(z,w) \to (\xi,\eta)$,
 we conclude that
 $$
 k^2 < 1.
 $$
 By
 $$
 \frac{d\xi}{dz} = \frac{\nu-\mu}{\lambda(z-\mu)^2},
 $$
 we have
 \begin{equation}
 \label{constC}
 \begin{split}
 \int \frac{dz}{w} = C \int \frac{d\xi}{\eta}, \\
 \int \frac{z dz}{w} = C \left(\mu \int \frac{d\xi}{\eta} + \frac{\mu-\nu}{\lambda} \int \frac{d\xi}{\left(\xi - \frac{1}{\lambda}\right)\eta}\right), \\
  C = \frac{\lambda(\nu-\mu)}{\sqrt{P(\nu)}}  >0.
 \end{split}
 \end{equation}

 Finally, let us introduce the Jacobi function $\sn(t,k)$:
 $$
u =  \int_0^\tau \frac{d\xi}{\sqrt{1-\xi^2}\sqrt{1-k^2\xi^2}} =
 \int_0^\theta \frac{d\varphi}{\sqrt{1-k^2 \sin^2\varphi}}
 $$
 where $\xi = \sin \varphi, \tau =  \sin \theta$,
 and we put
 $$
 \tau = \sn (u,k).
 $$
 This function is periodic:
 $$
 \sn (u+4K,k) = \sn (u,k)
 $$
 where
 $$
 K = \int_0^{\pi/2} \frac{d\varphi}{\sqrt{1-k^2 \sin^2\varphi}}.
 $$

\end{document}